\documentclass[11pt]{article}
\usepackage{amsfonts,epsf,amsmath,amssymb}
\usepackage{epsfig}
\usepackage{latexsym}
\usepackage{pgf,tikz}
\usepackage{bm}
\usepackage{lineno}
\usepackage{enumitem}

\newtheorem{theorem}{\bf Theorem}[section]
\newtheorem{corollary}[theorem]{\bf Corollary}
\newtheorem{lemma}[theorem]{\bf Lemma}
\newtheorem{proposition}[theorem]{\bf Proposition}

\newcommand{\proof}{\noindent{\bf Proof.\ }}
\newcommand{\qed}{\hfill $\square$ \bigskip}

\newcommand{\Fib}{{\cal F}}

\begin{document}
\title{A relation between Wiener index and Mostar index for daisy cubes}

\author{
Michel Mollard\footnote{Institut Fourier, CNRS, Universit\'e Grenoble Alpes, France email: michel.mollard@univ-grenoble-alpes.fr}
}
\date{\today}
\maketitle

\begin{abstract}
Daisy cubes are  a class of isometric subgraphs of the hypercubes $Q^n$. Daisy cubes include  some previously well known families of graphs like Fibonacci cubes and Lucas cubes. Moreover they appear in chemical graph theory.

Two distance invariants, Wiener and Mostar indices, have been introduced in the context of the mathematical chemistry. The Wiener index $W(G)$ is the sum of distance between  all unordered pairs of vertices of a graph $G$. The Mostar index $\mathit{Mo}(G)$ is a measure of how far $G$ is from being distance balanced.

In this paper we establish that the  Wiener and the Mostar indices of a daisy cube $G$ are linked by the relation
$ 2W(G)-\mathit{Mo}(G)=|V(G)||E(G)|. $ We deduce an expression of Wiener and Mostar index for daisy cubes. 
\end{abstract}

\noindent
{\bf Keywords:} daisy cube, Wiener index, Mostar index, Fibonacci cube, partial cube. 

\noindent
{\bf AMS Subj. Class. }: 05C07,05C35

\section{Introduction }
The {\em Fibonacci cube} of dimension $n$, denoted as $\Gamma_n$, is the subgraph of the hypercube $Q_n$  induced by vertices with no consecutive 1's. This graph was introduced in \cite{hsu93} as an interconnection network.

Lucas cubes~\cite{mupe2001}, are the cyclic version of Fibonacci cubes.
Structural properties of these graphs have been widely studied (see~\cite{survey} for a survey). 
Fibonacci cubes also play a role in mathematical chemistry. Indeed  they are precisely the resonance graphs of fibonacenes  an important class of hexagonal chains~\cite{KLZI-2005}. Later in~\cite{Z-2012}, a similar connection have been found  between Lucas cubes and the resonance graphs of  cyclic
polyphenanthrenes, which are related to non-cyclic fibonacenes.
Fibonacci cubes and Lucas cubes belong to daisy cubes~\cite{KM-2019} a familly of isometrical subgraphs of hypercubes.
The connection between daisy cubes and  resonance graphs of  catacondensed even ring systems have been explored in~\cite{BTZ-2020}.

The Wiener index $W(G)$ of a connected graph $G$ is the sum of distance all unordred pairs of vertices of $G$. This distance invariant is important in mathematical chemistry. The Wiener index of $\Gamma_n$ and $\Lambda_n$ have been determined in  ~\cite{KM-2012}.
Recently the Mostar index $\mathit{Mo}(G)$ has been introduced in  \cite{DMSTZ-2018} again in the context of graph chemical theory. It measures how far $G$ is from being distance-balanced.  The Mostar index of Fibonacci and Lucas cubes have been determined in \cite{ESS-2021}. In this note we prove that if $G$ is a daisy cube then the Wiener and the Mostar index are linked by the relation
$$ 2W(G)-\mathit{Mo}(G)=|V(G)||E(G)|.$$
  In the last section we derive similar expressions for $W(G)$ and $\mathit{Mo}(G)$ from the sequence of the number of edges using the direction  $i$ for $i\in[n]$.
\section{Preliminaries}
\label{sec:basic}
We will next give  some concepts and notations needed in this paper.
We denote by $[n]$ the set of integers $i$ such that $1\leq i \leq n$.
Let $\{F_n\}$ be the \emph{Fibonacci numbers}:
$F_0 = 0$, $F_1=1$, $F_{n} = F_{n-1} + F_{n-2}$ for $n\geq 2$. 
Let $B=\{0,1\}$. If will be convenient to identify elements $u = (u_1,\ldots, u_n)\in B^n$ and strings of length $n$ over $B$. We thus  briefly write $u$ as $u_1\ldots u_n$ and call $u_i$ the $i$th coordinate of $u$. We will use the power notation for the concatenation of bits, for instance $0^n = 0\ldots 0\in B^n$. 
We will denote by $\overline{u_i}$ the binary complement of $u_i$. 

The vertex set of $Q_n$, the \emph{hypercube of dimension $n$},  is the set $B^n$, two vertices being adjacent iff they differ in precisely one coordinate. We will say that an edge $uv$ of $Q_n$ uses the direction $i$ if $u$ and $v$ differ in the coordinate $i$. 

The \emph{distance} between two vertices $u$ and $v$ of  a graph $G$  is the 
number of edges on a shortest shortest $u,v$-path. It is immediate that the distance between two vertices of $Q_n$ is the number of coordinates the strings differ, sometime called Hamming distance.

The \emph{Wiener index} $W(G)$ of a connected graph $G$ is defined as the sum of all distances between pairs of vertices of $G$. Hence,
$$ W(G) =\sum_{\{u,v\} \subset V(G)}d(u,v).$$

A {\em Fibonacci string} of length $n$ is a binary string $b=b_1 b_2\ldots b_n$ with $b_i\cdot b_{i+1}=0$ for $1\leq i<n$. In other words a Fibonacci string is a binary string without $11$ as substring.\\ 
The {\em Fibonacci cube} $\Gamma_n$ ($n\geq 1$) is the subgraph of $Q_n$ induced by $\Fib_{n}$ the set of Fibonacci strings of length $n$. Because of the empty string $\epsilon$, $\Gamma_0 = K_1$. 

Not that for any integer $n$, $|V(\Gamma_{n})|=F_{n+2}$.

A Fibonacci string $b$ of length $n$ is a \emph{Lucas string} if $b_1 \, \cdotp b_n \neq 1$. 
That is, a Lucas string has no two consecutive 1s including the first and the last elements of the string. 
The \emph{Lucas cube} $\Lambda_n$ is the subgraph of $Q_n$ induced by the Lucas strings of length $n$. 
We have $\Lambda_0 = \Lambda_1 = K_1$.

Fibonacci cubes and Lucas cubes where extended to generalized Fibonacci cubes \cite{IKR-2012} and generalized Lucas cubes \cite{IKR-2012b}. For any arbitrary string $f$ the generalized Fibonacci cubes $Q_n[f]$ is the subgraph of $Q_n$induced by strings of  $B^n$ which do not contain $f$ as substring. Similarly the generalized Lucas cubes $Q_n[\stackrel{\leftarrow}{f}]$ is induced  by strings without a circulation containing $f$ as substring. Classical Fibonacci and Lucas cubes correspond to $f=11$.

If $u$ and $v$ are vertices of a graph $G$, the \emph{interval} $I_G(u,v)$ between $u$ and $v$ (in $G$) is the set of vertices lying on shortest $u,v$-path, that is, $I_G(u,v) = \{w | d(u,v) = d(u,w) + d(w,v)\}$. We will also write $I(u,v)$ when $G$ will be clear from the context. 
A subgraph $G$ of a graph $H$ is an \emph{isometric subgraph} 
if the distance between any vertices of $G$ equals the distance 
between the same vertices in $H$. 
Isometric subgraphs of hypercubes are called \emph{partial cubes}. 
The {\em dimension} of a partial cube $G$ is the smallest integer
$d$ such that $G$ is an isometric subgraph of $Q_d$. 
Many important classes of graphs are partial cubes, 
in particular trees, median graphs, benzenoid graphs, phenylenes, 
grid graphs and bipartite torus graphs. In addition, Fibonacci  
and Lucas cubes are partial cubes as well, see \cite{KL-2005}.

If $G=(V(G),E(G))$ is a graph and $X\subseteq V(G)$, then $\langle X\rangle$ denotes the subgraph of $G$ induced by $X$.  
Let $\le$ be a partial order on $B^n$ defined with $u_1\ldots u_n \le v_1\ldots v_n$ if $u_i\le v_i$ holds for $i\in [n]$. For $X \subseteq B^n$ we define the graph  $Q_n(X)$ as the subgraph of $Q_n$ with  
$$Q_n(X) = \left\langle \{u\in B^n| u\le x\ {\rm for\ some}\ x\in X \} \right\rangle$$
and say that $Q_n(X)$ is the {\em daisy cube generated by $X$}. Note that  if $\widehat{X}$ is the antichain consisting of the maximal elements of the poset $(X,\le)$, then $Q_n(\widehat{X}) = Q_n(X)$. As noticed in the daisy cube introductory paper \cite{KM-2019} we can alternatively say that $$Q_n(X) = \left\langle \bigcup_{x\in X} I_{Q_n}(x,0^n)\right\rangle=\left\langle \bigcup_{x\in \widehat{X}} I_{Q_n}(x,0^n)\right\rangle.$$
 Finally we will say that a graph $G$ is {\em a daisy cube} if there exist an isometrical embedding  of $G$ in some hypercube $Q_n$ and a subset $X$ of $B^n$ such that $G$ is the daisy cube generated by $X$. Such an embedding will be called a {\em proper embedding}.

By construction daisy cubes are partial cubes and as noticed in the same paper Fibonacci cubes, Lucas cubes, bipartite wheels, vertex-deleted cubes and hypercubes themselves are daisy cubes. It is easy to see that Pell graphs~\cite{M-2019} are also daisy cubes. Furthermore A. Vesel proved~ \cite{V-2019} that the cube complement of a daisy cube is a daisy cube.

For any fixed integer $s\geq 2$ the generalized Fibonacci cubes $Q_n[1^s]$ and Lucas cubes $Q_n[\stackrel{\leftarrow}{1^s}]$ are also daisy cubes. The Wiener index of these two families of graphs have been studied in \cite{KR-2015}.

A construction of daisy cubes in terms of expansion procedure is given in \cite{T-2020}. 

Let a partial cube $G$ of dimension $n$ be given together with its 
isometric embedding into $Q_n$. Then for $i=1,2,\ldots, n$ and 
$\chi =0,1$  the \emph{semicube} $W_{(i,\chi)}$ is defined as follows: 
$$W_{(i,\chi)}(G)=\{ u=u_1u_2\ldots u_n\in V(G)\ |\ u_i=\chi\}\,.$$
For a fixed $i$, the pair $W_{(i,0)}(G), W_{(i,1)}(G)$ of semicubes is 
called a \emph{complementary pair of semicubes}. The Wiener index of partial cubes can be determined using  the following result~\cite{KGM-1995}. 
\begin{theorem}\label{thm:basic}
Let $G$ be a partial cube of dimension $n$ isometrically 
embedded into $Q_n$. Then 
$$W(G) = \sum_{i=1}^n |W_{(i,0)}(G)|\cdot |W_{(i,1)}(G)|\,.$$
\end{theorem}
 Notice that, as expected, this expression is independent of the embedding. Indeed isometric embeddings of $G$ into $Q_n$ are 
unique up to the automorphisms of $Q_n$, see~\cite{WI-1984}. These automorphisms are generated by $\Theta=\{ \theta_{i,j} |i,j\in[n]\}$ and $T=\{ \tau_{i} |i\in [n]\} $ where  $\theta_{i,j}$ is the permutation of coordinates $i$ and $j$ and  $\tau_{i}$ is the complementation  of the coordinate $i$. A permutation  $\theta_{i,j}$ induces a similar permutation on complementary pairs and an element of $T$ is also without effect on the expression of $W(G)$ since $\tau_{i}$ exchanges $W_{(i,0)}(G)$ and $W_{(i,1)}(G)$.

Let $G=(V(G),E(G))$ be a connected graph. For any edge $uv\in E(G)$, let $n_{u,v}$ denote the number of vertices wich are closer to $u$ than to $v$. The \emph{Mostar index} $\mathit{Mo}(G)$  is defined as
$$ \mathit{Mo}(G) =\sum_{uv \in E(G)}|n_{u,v}-n_{v,u}|.$$
 An expression of the Mostar index for partial cubes in given in \cite{GRE-2020}.

\section{A relation}
This section is devoted to the proof of the following theorem.

\begin{theorem}
\label{thm:relation}
Let $G$ be a daisy cube. Then the Wiener and the Mostar indices of $G$ are linked by the relation
$$ 2W(G)-\mathit{Mo}(G)=|V(G)||E(G)|.$$
\end{theorem}

Like in~\cite{KM-2012} where it is applied in particular to Fibonacci cubes and Lucas cubes, theorem~\ref{thm:basic} will be our starting point for proving the relation.
Consider a proper embedding of $G$ in $Q_n$ and thus let  $X \subseteq B^n$ such that $G=\left\langle \{u\in B^n| u\le x\ {\rm for\ some}\ x\in X \} \right\rangle$. For any $i\in[n]$ let $E_i$ be the set of edges of $G$ using the direction $i$. 
 
\begin{proposition}
\label{le:Wgeq} For any $i\in[n]$ we have
$$|W_{(i,0)}(G)|\geq |W_{(i,1)}(G)|\,.$$

\end{proposition}
\proof
Let $u=u_1 u_2\ldots u_n$ in $W_{(i,1)}(G)$ and consider $\theta(u)=u_1\ldots u_{i-1} 0 u_{i+1}\ldots u_n$. Note that $\theta(u)\le u$  and since $u$ is a vertex of $G$ there exists $x\in X$ with $u \le x$. Therefore $\theta(u)\le x$ and $\theta(u)\in V(G)$. By this way we construct an injective mapping from $W_{(i,1)}(G)$ to $W_{(i,0)}(G).$
\qed

\begin{proposition}
\label{le:EieqW1} 
For any $i\in[n]$ we have $$|E_i|=|W_{(i,1)}(G)|\,.$$

\end{proposition}
\proof
Indeed let $u=u_1 u_2\ldots u_n$ in $W_{(i,1)}(G)$ and let $v=u_1\ldots u_{i-1} 0 u_{i+1}\ldots u_n$. It is clear that $v$ is a vertex of $G$ and that the edge $uv$ belongs to $E_i$. Reciprocally exactly one of the extremities of a given edge of $E_i$ belongs to $W_{(i,1)}(G)$. We obtain a bijective  mapping between $W_{(i,1)}(G)$ and $E_i$.
\qed
\begin{proposition}
\label{le:Eq2} 
For any edge $uv$ of $E_i$ with $u_i=0$ we have 
$$W_{(i,0)}(G)=\{w\in V(G) |d(w,u)<d(w,v)\}\,$$
$$W_{(i,1)}(G)=\{w\in V(G) |d(w,v)<d(w,u)\}\,.$$

\end{proposition}
\proof
Since $d(w,u)$ and $d(w,v)$ are the number of coordinates the strings differ and since $u,v$ differ only by the coordinate $i$ it is clear that $d(w,u)=d(w,v)+1$ if $w\in W_{(i,1)}(G)$ and $d(w,u)=d(w,v)-1$ otherwise.
\qed

\begin{lemma}\label{lm:MO}
Let $G$ be a daisy cube of dimension $n$ properly 
embedded into $Q_n$. Then 
$$\mathit{Mo}(G) = \sum_{i=1}^n |W_{(i,1)}(G)|(|W_{(i,0)}(G)|-|W_{(i,1)}(G)|)\,.$$
\end{lemma}
\proof
Let $e=uv$ be an edge of $E_i$ with $u_i=0$. By propositions \ref{le:Eq2} and \ref{le:Wgeq} we have $n_{u,v}=|W_{(i,0)}(G)| \geq |W_{(i,1)}(G)|=n_{v,u}.$
The contribution of the edge $e$ to $\sum_{uv \in E(G)}|n_{u,v}-n_{v,u}|$ is thus $|W_{(i,0)}(G)|-|W_{(i,1)}(G)|$. Therefore 
$$
\mathit{Mo}(G) =\sum_{i=1}^n\sum_{uv \in E_i}(|W_{(i,0)}(G)|-|W_{(i,1)}(G)|)
=\sum_{i=1}^n|E_i|(|W_{(i,0)}(G)|-|W_{(i,1)}(G)|).
$$
The conclusion follows from proposition \ref{le:EieqW1}.
\qed\\
We can now proceed to the proof of our relation. By theorem~\ref{thm:basic} we obtain $$W(G) = \sum_{i=1}^n |W_{(i,0)}(G)|\cdot |W_{(i,1)}(G)|\,.$$ Since $|W_{(i,0)}(G)|+|W_{(i,1)}(G)|=|V(G)|$ we deduce from Lemma~\ref{lm:MO} that
\begin{align*}
2W(G)-\mathit{Mo}(G) &= \sum_{i=1}^n |W_{(i,1)}(G)|(2|W_{(i,0)}(G)|-|W_{(i,0)}(G)|+|W_{(i,1)}(G)|)\\
&=\sum_{i=1}^n |W_{(i,1)}(G)|.|V(G)|\,.
\end{align*}
Since $E(G)$ is the disjoint union $E(G)=\cup_{i=1}^n E_i $ we have $\sum_{i=1}^n |W_{(i,1)}(G)|=|E(G)|$ and the relation follows.
\qed

\section{Conclusion}
The Wiener and Mostar indices of daisy cubes are completely determined by $|V(G)|$ and the sequence  $|E_i|$ (for $i\in[n])$ of the number of edges using the direction $i$  which is identical to the sequence of  $|W_{(i,1)}(G)|$. 
Indeed From theorem~\ref{thm:basic} and lemma~\ref{lm:MO} we have the relation $$W(G)-\mathit{Mo}(G) = \sum_{i=1}^n |W_{(i,1)}(G)|^2.$$
Combining this identity with that of theorem~\ref{thm:relation} we obtain the following assertion.

\begin{corollary}
\label{cor:indices}
Let $G$ be a daisy cube  properly embedded into $Q_n$. For $i\in[n]$ let $|E_i|$ be the number of edges using the direction $i$. Then the Wiener and the Mostar indices of $G$ are 
$$ W(G)=|V(G)||E(G)|-\sum_{i=1}^n |E_i|^2$$
$$ \mathit{Mo}(G)=|V(G)||E(G)|-2\sum_{i=1}^n |E_i|^2.$$
\end{corollary}

In conclusion of this paper note that it will be interesting to give bijective proofs of theorem~\ref{thm:relation} and its corollary.

\end{document}